\documentclass[12pt]{amsart}
\usepackage{amssymb,amsmath,amsthm,amscd}

\numberwithin{equation}{section}

\newtheoremstyle{fancy}{10pt}{10pt}{\itshape}{12pt}{\textsc\bgroup}{.\egroup}{8pt}{
}
\newtheoremstyle{fancy2}{10pt}{10pt}{}{12pt}{\itshape}{.}{8pt}{ }

\theoremstyle{fancy}

\newtheorem{lemma}[equation]{Lemma}
\newtheorem{proposition}[equation]{Proposition}
\newtheorem{theorem}[equation]{Theorem}

\newtheorem{main}{Theorem}

\theoremstyle{fancy2}

\newtheorem{remark}[equation]{Remark}

\newcommand{\cref}[1]{Corollary~\ref{#1}}

\newcommand{\taref}[1]{Table~\ref{#1}}

\newcommand{\no}{\noindent}

\newcommand{\C}{{\mathbb{C}}}

\newcommand{\Z}{{\mathbb{Z}}}
\newcommand{\Q}{{\mathbb{Q}}}

\newcommand{\SU}{\ensuremath{\operatorname{SU}}}

\newcommand{\diag}{\ensuremath{\operatorname{diag}}}

\newcommand{\Lk}{\ensuremath{\operatorname{Lk}}}

\def\con#1=#2(#3){#1 \equiv #2 \bmod{#3}}

\newcommand{\co}{{cohomogeneity }}
\begin{document}
\date{\today}

\title{Topological properties of Eschenburg spaces  and 3-Sasakian manifolds}

\author{Ted Chinburg}
\address{University of Pennsylvania\\
   Philadelphia, PA 19104}
\email{ted@math.upenn.edu}
\author{Christine Escher}
\address{Oregon State University\\
   Corvallis, OR 97331}
\email{tine@math.orst.edu}
\author{Wolfgang Ziller}
\address{University of Pennsylvania\\
   Philadelphia, PA 19104}
\email{wziller@math.upenn.edu}

\thanks{
The second author was supported by a grant from the Association for Women in
 Mathematics. The last author was supported by the Francis J. Carey Term Chair,
and the first and last author were supported by a grant from the National
Science Foundation.}

\maketitle

Riemannian manifolds with positive sectional curvature have been a
frequent topic of global Riemannian geometry for over 40 years.
Nevertheless, there are relatively few known examples of such
manifolds.  The purpose of this article is to study the
topological properties of some of these examples, the so-called
Eschenburg spaces, in detail.

In addition to positively curved metrics,  some Eschenburg spaces
also carry another special geometric structure, namely a 3-Sasakian
metric, i.e. a metric whose Euclidean cone is Hyperk\"ahler
\cite{BGMsurvey}.   3-Sasakian spaces are interesting since they are
Einstein manifolds and are connected to several other geometries:
They admit an almost free, isometric action by $SU(2)$ whose
quotient  is a quaternionic K\"{a}hler orbifold.  The twistor space
of this orbifold, which can be viewed as an $S^1$-quotient of the
3-Sasakian manifold, carries a natural K\"{a}hler-Einstein orbifold
metric with positive scalar curvature.

3-Sasakian structures are rare and  rigid, in fact the moduli space
of such metrics on a fixed manifold consists of at most isolated
points.  This motivated C. Boyer and  K. Galicki to pose the
question in \cite{BGMsurvey}[Question 9.9, p. 52] whether a manifold
can admit more than one 3-Sasakian structure.  Natural candidates
for such examples are the 3-Sasakian metrics discovered in
\cite{BGM}.  They are defined on the
 Eschenburg  biquotients $E_{a,b,c}=\diag(z^{a}, z^{b}, z^{c})\backslash
SU(3)/ \diag(z^{a+b+c}, 1,1)$, where $a,b,c$ are positive, pairwise
relatively prime integers.  The simplest topological invariant of
these spaces is the order of the fourth cohomology group, which is a
finite cyclic group of order $r=ab+ac+bc$.  By studying further
topological invariants of these manifolds we show:

\begin{main}
 For $r \le 10^7$, there  is a unique pair of $3$-Sasakian
Eschenburg spaces $E_{a,b,c}$ which are diffeomorphic to each other,
but not isometric.  It is given by $(a,b,c) = (2279, \,1603, \,384
)$ and $(2528 , \,939 , \,799 )$ with $r=5143925$.
\end{main}

The two $3$-Sasakian metrics are non-isometric since
  the isometric action by $SU(2)$ has cyclic
isotropy groups of order $a+b , a+c$ and $b+c$. In \cite{BGMsurvey}
they also asked whether two 3-Sasakian manifolds can be homeomorphic
to each other but not diffeomorphic. This happens frequently among
the 3-Sasakian Eschenburg spaces.
 There are 96 such pairs for
$r\le 10^7$, the first one of which is given by  $ (a,b,c) = (171,
\,164, \,1)$ and $(223, \,60, \,53 )$ for $r = 28379$. See
\taref{SH} for the next 4 such pairs.

The manifolds $E_{a,b,c}$ also carry a metric of positive sectional
curvature, although the 3-Sasakian metric never has positive
curvature.  They are special cases of the more general
 family  of Eschenburg spaces given by $E_{k,l}=\diag(z^{k_1}, z^{k_2},
z^{k_3})\backslash SU(3)/ \diag(z^{l_1}, z^{l_2}, z^{l_3})$.
 In this article we also examine the topology of the positively curved spaces among
this more general class of Eschenburg spaces.  They  contain in
particular the homogeneous Aloff-Wallach spaces
$SU(3)/\diag(z^p,z^q,\bar{z}^{p+q})$, \cite{AW}, which Kreck-Stolz
\cite{KSAloff-Wallach} classified up to homeomorphism and
diffeomorphism.  They were thus able to construct the first examples
of positively curved Riemannian manifolds which are homeomorphic but
not diffeomorphic. For this purpose they introduced three invariants
for $7$-dimensional manifolds with the same cohomology ring as
$W_{p,q}$ or $E_{k,l}$, which are generalizations of the classical
Eells-Kuiper invariant. They are computed using a bounding
$8$-dimensional manifold and detect both the homeomorphism and
diffeomorphism type.  In  case of the Aloff-Wallach spaces, which
can be viewed as
 circle bundles over the homogeneous flag
manifold, the bounding manifold is simply the corresponding disc
bundle. Another special case, namely the circle bundles over the
inhomogeneous $6$-dimensional flag manifold, were studied  in
\cite{AMP1},\cite{AMP2} although in this case there are not even any
homeomorphic pairs, see \cite{Sh}.  The invariants for the general
Eschenburg family, for which it is more difficult to find a bounding
$8$-manifold, were computed by Kruggel \cite{Kruggeldiffeo}.  We use
his formulas to study the topology of $E_{k,l}$.

The fourth cohomology group  of $E_{k,l}$ is  a finite cyclic group
of order
 $r=|k_1 k_2 +
k_1 k_3 + k_2 k_3 - (l_1 l_2 + l_1 l_3 + l_2 l_3)|$, and we
show that for a given value of r there are only finitely many
positively curved Eschenburg spaces.

\begin{main}
For $r \le 8000$, there is a unique pair of positively curved
Eschenburg spaces $E_{k,l}$ which are homeomorphic to each other,
but not diffeomorphic,  given by $(k_1,k_2,k_3 \;\vert\;
l_1,l_2,l_3) = (79, \,49, \,-50 \;\vert\; 0, \,46, \,32)$ and $(75,
\,54, \,-51 \;\vert\; 0, \,46 , \,32 )$ with $r=4001$.
\end{main}

 There are 69 pairs  of this type for $r\le 50000$, the first 5 are
listed in Table 4.2. Among these 69 there are also 4 pairs which are
diffeomorphic to each other, see Table 4.3. We were not able to show
that they are not isometric, although one should  expect this to be
the case. In the case of the Aloff-Wallach examples in
\cite{KSAloff-Wallach} the integer parameters must be significantly
larger. They find 11 homeomorphic pairs  and 3 diffeomorphic pairs
for $r<10^{17}$

To prove these theorems we use the Kreck-Stolz invariants, as
described in \cite{Kruggeldiffeo}.  In contrast to the case of the
circle bundles in \cite{KSAloff-Wallach} and \cite{AMP1}, where the
formulas for the invariants are fourth degree polynomials, the
formulas in \cite{Kruggeldiffeo} are quite complicated and involve
several number theoretic sums, with values in $\Q/\Z$. In order to
compute these invariants on a computer, one needs to control the
denominators, see Theorem 3.1.

It is a pleasure to thank D. Zagier and P.Gilkey for helpful
discussions and B. Bunker for assistance with the C-code.

\section{Eschenburg spaces}
\label{Eschenburg spaces}

A biquotient is a generalization of a homogeneous space where
$H\subset G\times G$ acts on $G$ via $(h_1,h_2)\cdot g =
h_1gh_2^{-1}$. The action is free if and only if $h_1$ is never
conjugate to $h_2$, in which case the quotient is a manifold denoted
by $G/\!/H$. We will also use the notation $\phi_1(h)\backslash G /
\phi_2(h)$ where the inclusion $H\subset G\times G$ is given by
$(\phi_1(h),\phi_2(h))$.

The Eschenburg spaces are an infinite family of seven dimensional manifolds
containing a subfamily that admits a metric a positive sectional curvature.
They were introduced by Eschenburg in 1982, \cite{E1}, and they can be described as
biquotients of $\SU(3)$.

Let $k:= (k_1,k_2,k_3)\,,\,l:= (l_1,l_2,l_3) \in \Z^3$ be two
triples of integers such that $k_1 + k_2 + k_3 = l_1 + l_2 + l_3$.
We can then define a two-sided action of $S^1 = \{z\in \C \mid
|z|=1 \}$ on $\SU(3)$  whose quotient we denote by $E_{k,l}$:
$$E_{k,l} :=  \diag(z^{k_1}, z^{k_2},
z^{k_3})\backslash SU(3)/ \diag(z^{l_1}, z^{l_2}, z^{l_3}) \,,\,
k_1 + k_2 + k_3 = l_1 + l_2 + l_3\,.$$ \no The action is free if and only if
$\diag(z^{k_1},z^{k_2},z^{k_3})$ is not
conjugate to $\diag(z^{l_1},z^{l_2},z^{l_3})$ which translates into the
following conditions, which must all be satisfied:
\begin{equation}
\label{freeness}
\begin{aligned}
\gcd(k_1-l_1,k_2-l_2) & = 1 \quad \gcd(k_1-l_2,k_2-l_1) &= 1 \\
\gcd(k_1-l_1,k_2-l_3) & = 1 \quad \gcd(k_1-l_2,k_2-l_3) &= 1 \\
\gcd(k_1-l_3,k_2-l_1) & = 1 \quad \gcd(k_1-l_3,k_2-l_2) &= 1
\end{aligned}
\end{equation}
In \cite{E1} it is also shown that the spaces $E_{k,l}$, equipped
with a  metric induced by a certain left invariant metric on
$SU(3)$, has positive sectional curvature if and only if:

\begin{equation}
\label{positive}
\begin{aligned}
k_i & \notin [\min(l_1,l_2,l_3),\max(l_1,l_2,l_3)]   &\text{for all} \quad i = 1, 2, 3\,.
\end{aligned}
\end{equation}

\begin{remark}  \label{group action}
One of the difficulties of dealing with these spaces, is that they
do not have a unique representation. One can easily change the
integers such that the $S^1$-actions are equivalent to each other,
and hence the quotient manifolds are diffeomorphic:

\begin{itemize} \item  We can use any permutation of the entries in
$k$ since an element of the Weyl group of $SU(3)$ acting on the
left will produce an equivalence of the corresponding actions. Similarly, we
 can use any permutation of the entries of $l$.

\item  We can switch all entries in $k$  with the entries in
$l$, if we replace the left-invariant metric with a right invariant one,
 since the inversion on $SU(3)$ induces an isometry.

\item  Simultaneously changing the signs of all entries in $k$
and
 $l$ is obtained by precomposing the action with $z\mapsto \bar z$.
Note though that in this case, the operation  changes the
orientation of $E_{k,l}$.

 \item  Adding an integer to all entries in $k$ and $l$,
 i.e. replacing  $k_i$ and $l_i$ by  $k_i + n$ and $l_i + n$ for $n \in
 \Z$, induces the same action of $S^1$, although in this case the
 new action may be only almost effective.
\end{itemize}
\end{remark}

In the case of positively curved Eschenburg spaces, we use
 such changes of the group action to obtain a unique
representation:

\begin{lemma} \label{representation}
Each positively curved Eschenburg space $E_{k,l}$ has the following unique representation
\begin{equation}
\begin{aligned}
k = &(k_1, k_2, l_1 + l_2 - k_1 - k_2)\\
l =  &(l_1, l_2, 0) \\
\text{with}  \quad &k_1 \ge k_2 > l_1 \ge l_2 \ge 0\,.
\end{aligned}
\end{equation}
\end{lemma}

\begin{proof}  Recall that the $k_i$ and $l_i$ have to satisfy the
positive curvature condition \eqref{positive}. If necessary change
the signs of all $k_i$ and $l_i$ to ensure that both  $k_1$ and
$k_2$ are on the right of the interval
$[\min(l_1,l_2,l_3),\max(l_1,l_2,l_3)]$.
 We can assume that $k_1 \ge
k_2$ by changing the order of $k_1$ and $k_2$ if necessary.  Now
subtract $\min(l_1,l_2,l_3)$ and after possibly changing the order
of $l_1$ and $l_2$ we can assume that $l_1 \ge l_2 \ge l_3 = 0$.
Hence we obtain that $k_1 \ge k_2 > l_1 \ge l_2 \ge l_3=0 > k_3 =
l_1 + l_2 - k_1 - k_2$.
\end{proof}

 Using the Serre spectral sequence
as in \cite{biquotients}, or by following the methods developed in
\cite{singhof}, one obtains the cohomology ring with integer
coefficients:

\begin{equation}
\label{cohomology}
\begin{aligned}
H^1(E_{k,l}) & = 0 \,, \,   H^2(E_{k,l})   = \Z  \quad  \text{generated by} \quad u \,,\\
H^3(E_{k,l}) & = 0 \,, \,   H^4(E_{k,l})   = \Z_{|r|} \quad  \text{generated by} \quad u^2 \,,\\
\text{with} \quad r & = \sigma_2(k) -  \sigma_2(l)\,,
\end{aligned}
\end{equation}
where $\sigma_i(k):=\sigma_i(k_1,k_2,k_3)$ is the $i$-th
elementary symmetric function.  Moreover, by studying the cell
structure of the Eschenburg spaces as in \cite[Remark
1.3]{homotopy}, one proves that $r$ is always an odd number.

In order to generate a complete list of all positively curved
Eschenburg spaces the following proposition will be important:

\begin{proposition} \label{finite}
For each odd $r \in \Z$ there are only finitely many positively
curved Eschenburg spaces $E_{k,l}$ with $H^4(E_{k,l}) = \Z_{|r|}$.
\end{proposition}

\begin{proof}
Assume that $E_{k,l}$ is represented as in Lemma
\eqref{representation}.  Then we obtain that
$$\begin{aligned}  r & = k_1\,k_2 + k_1 \, k_3 + k_2 \, k_3 - l_1 \, l_2  \quad &&\text{with} \quad  k_3 = l_1 + l_2 - k_1 - k_2\\
& = - [k_1 \, (k_1 - l_1) + (k_2 - l_2) \, (k_1 + k_2 - l_1)] \quad  &&\text{with} \quad k_1 \ge k_2 > l_1 \ge l_2 \ge 0\,.
\end{aligned}$$
Note that $k_1 > 0\,,\,k_1 - l_1 > 0\,,\, k_2 - l_2 > 0$ and $k_1
+ k_2 - l_1 > 0$ and hence $r < 0$. If we fix a positive odd
integer $N \in \Z^{+}$ with $N =  -[k_1 \, (k_1 - l_1) + (k_2 -
l_2) \, (k_1 + k_2 - l_1)]$, then the above conditions imply
 $k_1\,,\,k_2 \,,\, l_1\,,\,l_2 \in [1, N]$.  Hence there are only finitely many
choices for $k_i\,,\,l_i$ as claimed.
\end{proof}

\medskip

There are various interesting subfamilies of the Eschenburg spaces that have
 appeared in
  different contexts, see for example
 \cite{KSAloff-Wallach, AMP1,AMP2, BGM,D, Sh}. We use \cite{Z}
 for a systematic description of these subfamilies.

\medskip

\no (1)  Cohomogeneity one Eschenburg spaces

\smallskip

A group action of $G$ on a manifold $M$   is said to be of
cohomogeneity one if  the orbit space $M/G$ is one-dimensional.
For $E_{k,l}$ with $k_1 = k_2$ and $l_1 = l_2$ the group
$SU(2)\times SU(2)$ acting on $SU(3)$ on the left and on the
right, clearly commutes with the $S^1$-action and induces a \co
one action on $E_{k,l}$. Using a change in the group action as in
 \eqref{group action}, we can rewrite these cohomogeneity one
 Eschenburg spaces as follows.
$$ E_a= \diag(z^a,z,z) \backslash
\SU(3) /  \diag(z^{a+2},1,1)  \,.$$ The action is free for all
$a\in \Z $. Since $E_{a}=E_{-a-1}$, again via changes as in
\eqref{group action}, we can assume $a\ge 0$. It follows that
$E_{a}$ has positive curvature for all $a \ne 0$.
 Note that in
this case $r = 2\,a +1$ and we
obtain exactly one positively curved cohomogeneity one space for
each odd $r$.

\medskip

\no  (2) Cohomogeneity two Eschenburg spaces

\smallskip

If  two of the integers in $k$ or in $l$ are equal, we obtain an
action of $G=SU(2)\times T^2$ on  $SU(3)$ commuting with the $S^1$
action such that  the orbit space $E_{k,l}/G$ is two dimensional.
There are two families of Eschenburg spaces of this type.

$(2^+)$ $k_1 = k_2\,.$  We can rewrite these particular cohomogeneity
two  Eschenburg spaces as follows.
$$ E_{a,b,c}=\diag(z^a,z^b,z^c) \backslash
\SU(3) /  \diag(z^{a+b+c},1,1)$$ The action is free if $a,b,c$ are
pairwise relatively prime integers and $ E_{a,b,c}$ has positive
curvature if and only if
 $a \ge b \ge c > 0$.  We can assume
 $a>b>c>0$ since it is otherwise an Eschenburg space of cohomogeneity one. This
subfamily also includes the circle bundles over the inhomogeneous
flag manifold
 given by $a=b+c$, see \cite{E4}, and their topological properties were studied in
 \cite{AMP1,AMP2, Sh}.
  As mentioned in the introduction, they
   also admit a second metric which is $3$-Sasakian and we will
   study this subfamily to prove Theorem A.
\smallskip

$(2^-)$  $l_1 = l_2\,.$  These spaces can be rewritten again  as $$
\diag(z^a,z^b,z^c) \backslash \SU(3) /  \diag(z^{a+b+c},1,1)$$ with
$a,b,c$ pairwise relatively prime. But in this case they have
positive curvature if and only if
  $a \ge b >0 $ and $  -b < c < 0 $, and $a>b$ if they are not of cohomogeneity one.
If $a+b+c=0$ we recover the  Aloff-Wallach spaces $W_{a,b}
=SU(3)/\diag(z^a,z^b,z^{-a-b})$.

\medskip

\no  (3) Cohomogeneity four Eschenburg spaces

\smallskip
In the general case $G=T^4$ and  $M/G$ is four dimensional. In our
normalization  we can assume  $k_1
> k_2 > l_1 > l_2 \ge 0\,.$

In \cite{GSZ} it is shown that in the cohomogeneity one and two case,
the above group actions cannot be
extended to an action with smaller cohomogeneity.

\section{Topological invariants}
\label{Invariants}

The topological invariants we use are the order of the fourth
cohomology group $r:=|r(k,l)|= |\sigma_2(k) - \sigma_2(l)| $ , the
self-linking number, the first Pontrjagin class, and the
Kreck-Stolz invariants. All of these invariants were computed for
most of the Eschenburg spaces in \cite{Kruggeldiffeo}.  The
results can be summarized as follows.

\smallskip

(1) The self-linking number of a class $u^2 \in H^4(E_{k,l})$ is given by
 $$\Lk(k,l):= \Lk(E_{k,l}) = \Lk(u^2, u^2) = - \frac{s^{-1}(k,l)}{r(k,l)} \in
 \Q / \Z$$
 where $s(k,l):=\sigma_3(k) - \sigma_3(l)$
 and $s^{-1}(k,l)$ is the multiplicative inverse of $s(k,l)$ in $\Z_{|r(k,l)|}$.
 Notice that $\Lk(k,l)$ is uniquely determined by $s(k,l)$ up to
 sign (which will effect the orientation). In our normalization we
 always have $r(k,l) < 0$ and hence in our tables we will use $s(k,l) \mod |r(k,l)|$.

\smallskip

(2)  The first Pontrjagin class, as an element in $H^4(E_{k,l}) \cong  \Z_{|r(k,l)|}\,,$ is
given by
$$
p_1(k,l):= p_1(E_{k,l})  = [2\,\sigma_1(k)^2 - 6\,\sigma_2(k)]\cdot u^2 \, \in \Z_{|r(k,l)|}\,.$$

Note that the roles of $k$ and $l$ are interchangeable since $\sigma_1(k) = \sigma_1(l)$ and
$\sigma_2(k) \equiv  \sigma_2(l) \mod |r(k,l)|$.  The second Stiefel-Whitney class
vanishes for all $E_{k,l}$.

\smallskip

(3)  The Kreck-Stolz invariants are based on the Eells-Kuiper $\mu$-invariant and
are defined as linear combinations of relative characteristic numbers of appropriate
bounding manifolds.  They were introduced in \cite{KSAloff-Wallach} and calculated for
 most of the Eschenburg spaces in \cite{Kruggeldiffeo}.  One first constructs a cobordism of
$SU(3)$ and  extends the $S^1$-action
  to this cobordism.  However, in general this extended circle action
 is not free anymore and in order to make the action almost free, Kruggel introduced the following
  condition (C).

  We say that condition (C) holds
  if and only if the matrix
$A_{i,j}=(k_i - l_j)$ contains at least one row or one column
whose entries are pairwise relatively prime. In
\cite{Kruggeldiffeo}  it was indicated that condition (C) might
always be satisfied.    Unfortunately, this is not the case. There
are many Eschenburg spaces, even positively curved ones, which do
not satisfy condition (C).  For example for $r<5000$ there are $54$ positively curved
Eschenburg spaces for which condition (C) fails.

  The positively curved Eschenburg space
with smallest $|r(k,l)|$ where this occurs is  given by $
(k_1,k_2,k_3 \;\vert\; l_1,l_2,l_3)=  (35,  21,  -34, \;\vert\; 12,
10 ,  0 )$, with:

$$
 A=\begin{pmatrix} 23&  5^2  &  5 \cdot 7  \\
 3^2 &  11  &  3 \cdot 7  \\
 - 2 \cdot 23 &  - 2^2 \cdot 11 &  - 2^2 \cdot 17
 \end{pmatrix}
$$

Determining the homeomorphism type and diffeomorphism type of
Eschenburg spaces that do not satisfy condition (C) remains an
open problem.
 Notice though that in the subclass of cohomogeneity two Eschenburg spaces,
 and in particular
 for $3$-Sasakian manifolds,
condition (C)  always holds since in that case the first column $(a,
b, c)$ consists of pairwise relatively prime integers. Fortunately,
in   the proof of Theorems B this problem
does not arise either, as explained in section 4.

Assuming that condition (C) holds for the $j$-th column, there are at
most three exceptional orbits for the $S^1$ action on the cobordism of $SU(3)$
 with isotropy groups $\Z_{|k_1-l_j|}$,
$\Z_{|k_2-l_j|}$ and $\Z_{|k_3-l_j|}\,.$  After removing small
equivariant neighborhoods of these orbits the action becomes free
and the quotient is
 a smooth eight dimensional manifold $W_{k,l}$ with boundary
$\partial(W_{k,l}) = E_{k,l} \cup L_1 \cup L_2 \cup L_3$ where the $L_i$ are the
following seven dimensional lens spaces.

$$\begin{aligned}
L_1 & := L(k_1 - l_j; k_2 - l_j, k_3 - l_j, k_2 - l_{[j+1]_2}, k_3 - l_{[j+1]_2}) \\
L_2 & := L(k_2 - l_j; k_1 - l_j, k_3 - l_j, k_1 - l_{[j+1]_2}, k_3 - l_{[j+1]_2}) \\
L_3 & := L(k_3 - l_j; k_1 - l_j, k_2 - l_j, k_1 - l_{[j+1]_2}, k_2 - l_{[j+1]_2}) \\
\end{aligned}
$$
where we used the notation $[n]_p := m$ if $n = \lambda \cdot p + m$
for $m = 1, \dots, p$ , for the residue class $[n]$ modulo $p$.

Since the invariants are additive with respect to unions, we obtain
$s_i(W_{k,l}) = s_i(E_{k,l}) + s_i(L_1) + s_i(L_2) + s_i(L_3) \in \Q /
\Z\,,\,i= 1,2\,.$  Calculating $s_i(W_{k,l})$,
 yields the following expressions for the Kreck-Stolz invariants, which  hold in the case  condition (C) is satisfied for the $j$-th column.

\begin{equation}
\label{invariants}
 \begin{aligned}
s_{1}(E_{k,l}) &=  \frac{4\,|\,r(k,l)\,(k_1 - l_j)\,(k_2 -
l_j)\,(k_3 - l_j)\,| - q(k,l)^2}{2^{7}  \cdot 7  \cdot r(k,l)
\,(k_1 - l_j)\,(k_2 - l_j)\,(k_3 - l_j)}   \\
& - \sum_{i=1}^3 \, s_1(k_{[i]_3} - l_j; k_{[i+1]_3} - l_j, k_{[i+2]_3} - l_j, k_{[i+1]_3} - l_{[j+1]_2}, k_{[i+2]_3} - l_{[j+1]_2}) \\
s_{2}(E_{k,l}) &=  \frac{q(k,l) - 2}{2^{4} \cdot 3  \cdot  r(k,l) \,(k_1 - l_j)\,(k_2 - l_j)\,(k_3 - l_j)}  \\
& - \sum_{i=1}^3 \, s_2(k_{[i]_3} - l_j; k_{[i+1]_3} - l_j, k_{[i+2]_3} - l_j,
k_{[i+1]_3} - l_{[j+1]_2}, k_{[i+2]_3} - l_{[j+1]_2}) \end{aligned}
\end{equation}

\no where
$q(k,l):= (k_1 - l_j)^2 + (k_2 - l_j)^2 + (k_3 - l_j)^2+  (k_1 - l_{[j+1]_2})^2 + (k_2 - l_{[j+1]_2})^2 + (k_3 -
                l_{[j+1]_2})^2 -(l_j -l_{[j+1]_2})^2\,$, and
$s_i(p;p_{1},p_{2},p_{3},p_{4}):=   s_{i}(L_{p}(p_{1},p_{2},p_{3},p_{4}))
\in \Q / \Z\,,\,i=1,2$
are the Kreck-Stolz invariants of the lens space
$L_{p}(p_{1},p_{2},p_{3},p_{4})=S^7/\Z_p\,.$

The freeness condition of the action \eqref{freeness} implies that
$k_i - l_j \ne 0$ for $i,j=1,2,3$, and hence the above expressions for $s_1$ and
$s_2$ are well-defined.

By the Atiyah-Patodi-Singer index theorem the Kreck-Stolz
invariants can also be expressed as linear combinations of
eta-invariants.  Calculating these eta-invariants for the lens
spaces, one obtains:

\begin{alignat}{2} \label{lens space}
s_1(p;p_{1},p_{2},p_{3},p_{4}) &=\frac{1}{2^{5} \cdot 7 \cdot p}\,
\sum_{k=1}^{|p|-1}\,\prod_{j=1}^{4}\, \cot\left(\frac{k \pi
p_{j}}{p}\right) +  \frac{1}{2^{4} \cdot p} \sum_{k=1}^{|p|-1}\,
\prod_{j=1}^{4}\,\csc\left(\frac{k \pi p_{j}}{p}\right)\,;
\\
s_2(p;p_{1},p_{2},p_{3},p_{4}) &=  \frac{1}{2^{4} \cdot
p}\,\sum_{k=1}^{|p|-1}\, (e^{\frac{2 \pi \imath k}{|p|}} -1)\,
\prod_{j=1}^{4}\,\csc\left(\frac{k \pi p_{j}}{p}\right)\,. \notag \\
\notag \end{alignat}

\no For $p=\pm 1$ these expressions are interpreted to be $0$.

These formulas only hold in the case that condition (C) is
satisfied for the jth column.  We also list the case
where the jth row  consists of relatively prime entries, since this
will be needed in our calculations.

\begin{equation}
 \begin{aligned}
s_{1}(E_{k,l}) &=  \frac{4\,|\,r(k,l)\,(k_j - l_1)\,(k_j -
l_2)\,(k_j - l_3)\,| - q(k,l)^2}{2^{7}  \cdot 7  \cdot r(k,l)
\,(k_j - l_1)\,(k_j - l_2)\,(k_j - l_3)}   \notag\\
& + \sum_{i=1}^3 \, s_1(k_j - l_{[i]_3}; k_j - l_{[i+1]_3}, k_j - l_{[i+2]_3}, k_{[j+1]_2} - l_{[i+1]_3}, k_{[j+1]_2} - l_{[i+2]_3}) \\
s_{2}(E_{k,l}) &=  \frac{q(k,l) - 2}{2^{4} \cdot 3  \cdot  r(k,l) \,(k_j - l_1)\,(k_j - l_2)\,(k_j - l_3)}  \notag\\
& + \sum_{i=1}^3 \, s_2(k_j - l_{[i]_3}; k_j - l_{[i+1]_3}, k_j - l_{[i+2]_3}, k_{[j+1]_2} - l_{[i+1]_3}, k_{[j+1]_2} - l_{[i+2]_3})
\notag
\end{aligned}
\end{equation}
where $q(k,l):= (k_j - l_1)^2 + (k_j - l_2)^2 + (k_j - l_3)^2+  (k_{[j+1]_2} - l_1)^2 + (k_{[j+1]_2} - l_2)^2 + (k_{[j+1]_2} -  l_3)^2 -(k_j -k_{[j+1]_2})^2\,.$

Using these invariants we now state the classification theorems by Kruggel for the Eschenburg spaces.

\begin{theorem} (Kruggel) \cite{homotopy, Kruggeldiffeo} \label{classification}
Assume the Eschenburg spaces $E_{k,l}$ and $E_{k',l'}$ both satisfy condition (C).  Then

(I)  $E_{k,l}$ and $E_{k',l'}$ are (orientation preserving)
homeomorphic if and only if
\begin{alignat}{2}
& (a)  \hspace{10pt}  | r(k,l)| &&= |r(k',l')| \in \Z \notag \\
&(b)     \hspace{10pt} \Lk(k,l)  && \equiv \Lk(k',l') \in \Q / \Z \notag\\
&(c)    \hspace{10pt} p_1(k,l)&& \equiv p_1(k',l')\in \Z_{|r(k,l)|}  \notag \\
&(d)    \hspace{10pt}  s_2(E_{k,l})  && \equiv  s_2(E_{k',l'})\in \Q / \Z
\notag
\end{alignat}

(II) $E_{k,l}$ and $E_{k',l'}$ are (orientation preserving)
diffeomorphic if and only if  in addition
$$s_1(E_{k,l}) \equiv s_1(E_{k',l'})\in \Q / \Z$$

(III) $E_{k,l}$ and $E_{k',l'}$ are (orientation preserving)
homotopy equivalent if and only if
\begin{alignat}{2}
& (a) \hspace{50pt} | r(k,l)| &&= |r(k',l')| \in \Z \notag \\
&(b)\hspace{50pt} \Lk(k,l)&& \equiv \Lk(k',l) \in \Q / \Z   \notag\\
&(c) \hspace{50pt} s_{22}(E_{k,l})&& \equiv  s_{22}(E_{k',l'})\in \Q / \Z
\notag
\end{alignat}
where $ s_{22}(E_{k,l}):= 2\,|r(k,l)| \, s_2(E_{k,l})\,.$
\end{theorem}

\smallskip

For the corresponding theorem in the orientation reversing case the
linking number and the Kreck-Stolz invariants change signs. Recall
that in this theorem, $r(k,l)= \sigma_2(k) - \sigma_2(l) $,
$p_1(E_{k,l})  = [2\,\sigma_1(k)^2 - 6\,\sigma_2(k)]\cdot u^2 \, \in
\Z_{|r|}\,$, and the equality of the linking forms can replaced by
the equality of the numbers  $s(k,l):=\sigma_3(k) - \sigma_3(l) \in
\Z_{|r(k,l)|}$.

\begin{remark}
In \cite{KSAloff-Wallach} Kreck-Stolz used another invariant $s_3$
in the homeomorphism classification, and showed that $r(k,l) , s_2$
and $s_3$ determine the homeomorphism type.
Following \cite{Kruggeldiffeo}, the formula for the invariant $s_3$
for the Eschenburg spaces is easily seen to be:

$$\begin{aligned}
s_{3}(E_{k,l}) &=  \frac{q(k,l) - 8}{2^2\cdot 3  \cdot  r(k,l)
\,(k_1 - l_j)\,
(k_2 - l_j)\,(k_3 - l_j)}  \\
 &- \sum_{i=1}^3 \, s_3(k_{[i]_3} - l_j; k_{[i+1]_3} - l_j,
k_{[i+2]_3} - l_j, k_{[i+1]_3} - l_{[j+1]_2}, k_{[i+2]_3} -
l_{[j+1]_2}) \end{aligned}$$ where
$$s_3(p;p_{1},p_{2},p_{3},p_{4}) =  \frac{1}{2^{4} \cdot
p}\,\sum_{k=1}^{|p|-1}\, (e^{\frac{4 \pi \imath k}{|p|}} -1)\,
\prod_{j=1}^{4}\,\csc\left(\frac{k \pi p_{j}}{p}\right)$$

\end{remark}

\smallskip

\section{Number Theory}
\label{Number Theory}

The difficulty in explicitly calculating the invariants comes from
the complicated expressions  \eqref{lens space} for the Kreck-Stolz
invariants of the lens spaces.  These expressions are rational since
they lie in $\Q / \Z$. However, in order to use a computer program
to calculate the invariants it is  necessary to control the
denominators.    The bounds on the denominators in the theorem below
are similar to those obtained by Zagier in \cite{Zagier} for the
 higher order Dedekind sum $T$.  However, for the other
sums the results we need are not contained in \cite{Zagier}, so we
give a proof.

\begin{theorem}
\label{thm:denomtheorem} For all integers $p, p_1, \cdots , p_4$ such that $p$ is relatively prime to each $p_i$,
 the numbers
$$
T = \sum_{k=1}^{|p|-1}\, \prod_{j=1}^{4}\, \cot\left(\frac{k
\pi p_{j}}{p}\right) \hspace{5pt} , \hspace{5pt}  S = \sum_{k=1}^{|p|-1}\,
\prod_{j=1}^{4}\,\csc\left (\frac{k
\pi p_{j}}{p}\right )$$

$$R =  \sum_{k=1}^{|p|-1}\, \cos\left(\frac{2 \pi k
}{|p|}\right)\,
\prod_{j=1}^{4}\,\csc\left (\frac{k \pi p_{j}}{p}\right ) \hspace{5pt} \mathrm{and} \hspace{5pt}
U =  \sum_{k=1}^{|p|-1}\, \cos\left(\frac{4 \pi k
}{|p|}\right)\,
\prod_{j=1}^{4}\,\csc\left (\frac{k \pi p_{j}}{p}\right )
$$
are rational with denominators which divide 45. If $\sum_{i = 1}^4 p_i$ is
even,

$$
 \sum_{k=1}^{|p|-1}\, (e^{\frac{2 \pi \imath k}{|p|}} -1)\,
\prod_{j=1}^{4}\,\csc\left (\frac{k \pi p_{j}}{p}\right ) = R - S
\text{ and }
 \sum_{k=1}^{|p|-1}\, (e^{\frac{4 \pi \imath k}{|p|}} -1)\,
\prod_{j=1}^{4}\,\csc\left (\frac{k \pi p_{j}}{p}\right ) = U - S.
$$
\end{theorem}

\smallskip 

To begin the proof, note that the last two equalities are clear from
grouping terms for $k$ and $|p|-k$. For the rest of the proof, we
can assume that $p > 0$ after replacing $p$ by $-p$ if necessary.

Define

\begin{equation}\label{eq:Fdef}
\begin{aligned}
F_1(x) &= \prod_{j = 1}^4 \mathrm{cot}( x p_j\pi )\quad , \quad
&&F_2(x) = \prod_{j = 1}^4
\mathrm{csc}( xp_j\pi)\\
F_3(x) &= \mathrm{cos}(2 x \pi ) F_2(x) \quad , \quad &&F_4(x) =
\mathrm{cos}(4 x \pi ) F_2(x).
\end{aligned}
\end{equation}

If $F$ is one of the  $F_i$ then $F(x) = F(2 - x)$.
It follows that
\begin{equation}
\label{eq:Fzap}
\sum_{k = 1, k \ne p}^{2p-1} F(k/p) = 2T, \ 2S, \ 2R, \ 2U \quad {\rm if}\quad F = F_1, \ F_2, \ F_3, \ F_4 \quad {\rm respectively.}
\end{equation}

 For each divisor $p'$ of $2p$ define
\begin{equation}\label{eq:divnice}
\begin{aligned}
\mathcal{C}_{p'} &= \{k: 1 \le k < 2p  \quad \mathrm{and} \quad
(k,2p) = p'\}
\\
\mathcal{D}_{2p/p'} &= \{k': 1 \le k' < \frac{2p}{p'} \quad
\mathrm{and} \quad (k',\frac{2p}{p'}) = 1\}
\end{aligned}
\end{equation}
The map $k' \to k' \cdot p'$ defines a bijection between $\mathcal{D}_{2p/p'}$ and
$\mathcal{C}_{p'}$.  Therefore for $F \in \{F_1,F_2,F_3, F_4\}$,
\begin{eqnarray}
\label{eq:computationT}
\sum_{k = 1, k \ne p}^{2p-1} F(k/p) &=&
\sum_{p \ne p' | 2p} \ \sum_{k \in \mathcal{C}_{p'}} \
 F(k/p) \nonumber\\ &=&
\sum_{2 \ne d = \frac{2p}{p'}| 2p} \ \sum_{k' \in \mathcal{D}_{d}} \
F(2k'p'/dp')\\
 &=&
\sum_{2 < d = \frac{2p}{p'}| 2p} \ \sum_{k' \in \mathcal{D}_{d}} \
F(2k'/d)\nonumber
\end{eqnarray}
In deducing the last line, we use the fact
that $\mathcal{D}_{1}$ is empty.

Let $\zeta_{d} = \mathrm{exp}\left (\frac{2\pi i}{d}\right )$ be a
$d^{th}$ root of unity. For each divisor  $d$ of $2p$ with $d>2$,
denoted by  $2< d|2p  $,  there is a bijection $\mathcal{D}_d \to
\mathrm{Gal}({\Q}(\zeta_d)/{\Q})= G_d$ which sends $k' \in
\mathcal{D}_d$ to the unique automorphism $\sigma_{k'}$ of
${\Q}(\zeta_d)$ over ${\Q}$ which takes $\zeta_d$ to $\zeta_d^{k'}$.
From (\ref{eq:Fdef})  we have

\begin{equation}
\label{eq:Fdefeasy}
 F_1(2k'/d) = \prod_{j = 1}^4
\frac{\zeta_d^{k'p_j} + \zeta_d^{k'p_j}}{\zeta_d^{k'p_j} -
\zeta_d^{k'p_j}} \quad ; \quad F_2(2k'/d) = 16 \prod_{j = 1}^4
\frac{1}{\zeta_d^{k'p_j} - \zeta_d^{k'p_j}} \quad
\end{equation}
\begin{equation*}\label{eq:Fdefeasy2}
 F_3(2k'/d) =
\left ( \frac{\zeta_d^{2k'p_j} + \zeta_d^{2k'p_j}}{2}\right )
F_2(2k'/d)\quad ; \quad F_4(2k'/d) = \left ( \frac{\zeta_d^{4k'p_j}
+ \zeta_d^{4k'p_j}}{2}\right ) F_2(2k'/d).
\end{equation*}

Therefore for $F \in \{F_1,F_2,F_3, F_4\}$ and $k' \in
\mathcal{D}_d$ we have $F(2k'/d) \in {\Q}(\zeta_d)$ and
$\sigma_{k'}(F(2/d)) = F(2k'/d)$.  Hence
\begin{equation}
\label{eq:computationTnew}
\sum_{k = 1, k \ne p}^{2p-1} F(k/p) =
\sum_{2 < d = \frac{2p}{p'}| 2p} \ \sum_{k' \in \mathcal{D}_{d}} \
F(2k'/d)
=
\sum_{2 < d = \frac{2p}{p'}| 2p} \ \mathrm{Tr_d}(F(2/d))
\end{equation}
where $\mathrm{Tr_d}:{\Q}(\zeta_d) \to {\Q}$ is the trace function
defined by $\mathrm{Tr_d}(\tau) = \sum_{\sigma \in G_d}
\sigma(\tau)$.

Define
\begin{equation}
\label{eq:greekdef}
\alpha_d = \prod_{j = 1}^4  \left (
\zeta_{d}^{p_j} -
\zeta_{d}^{-p_j} \right ) \  , \  \beta_d = \prod_{j = 1}^4  \left (
\zeta_{d}^{p_j} +
\zeta_{d}^{-p_j} \right )\ , \ \gamma_d =  \zeta^2_d + \zeta^{-2}_d
, \ \tau_d =  \zeta^4_d + \zeta^{-4}_d.
\end{equation}
Putting together (\ref{eq:Fzap}), (\ref{eq:computationTnew}) and
(\ref{eq:Fdefeasy})
 shows

\begin{equation}
\label{eq:summary}
\begin{aligned}
2 T  &=  \sum_{2 < d | 2p} \mathrm{Tr_d} \left
(\frac{\beta_d}{\alpha_d} \right )\quad , \quad &&2S = 16 \sum_{2 <
d | 2p}  \mathrm{Tr_d} \left (\frac{1}{\alpha_d} \right )\ , \
\\
 2R &=  8 \sum_{2 < d | 2p} \mathrm{Tr_d}\left
 (\frac{\gamma_d}{\alpha_d} \right ) \quad , \quad &&2U =  8 \sum_{2 < d | 2p} \mathrm{Tr_d}
\left (\frac{\tau_d}{\alpha_d} \right ).
\end{aligned}
\end{equation}

\begin{lemma}
\label{lem:dump} Fix $2 < d | 2p$ and let $d' = d/(d,2) > 1$. Then
each of $\alpha_d$, $\beta_d$ , $\gamma_d$ and $\tau_d$ are integers
in the
 real subfield ${\Q}(\zeta_{d'})^+ = {\Q}(\zeta_{d'} + \zeta_{d'}^{-1})$
of ${\Q}(\zeta_{d'})$. Moreover,
\begin{enumerate}
\item{}
If $d'$ is not a prime power, then $\alpha_d$ is a unit.
\item{} Suppose $d' = l^r$ for some prime $l$ and some $r > 0$.  Let $\mathcal{P}$
be the unique prime ideal over $l$ in ${\Q}(\zeta_{d'})^+$.
 If $d' \ne 2$ then $\alpha_d$ generates the ideal $\mathcal{P}^2$, while if $d' = 2$ then
$\alpha_d$ generates $\mathcal{P}^4$. The relative degree
\hbox{$[{\Q}(\zeta_d):{\Q}(\zeta_{d'})^+]$} is $2$ if $l > 2$ or $d'
= 2$, and \hbox{$[{\Q}(\zeta_d):{\Q}(\zeta_{d'})^+] = 4$} if $l = 2$
and $d' > 2$.
\end{enumerate}
\end{lemma}

\begin{proof} Let us first show $\alpha_d$, $\beta_d$,
$\gamma_d$ and $\tau_d$ are integers in ${\Q}(\zeta_{d'})^+$. Since
$\zeta_d$ is integral, so are $\alpha_d$, $\beta_d$, $\gamma_d$ and
$\tau_d$.  Complex conjugation acts on ${\Q}(\zeta_d)$ by sending
$\zeta_d$ to $\zeta_d^{-1}$.  Hence (\ref{eq:greekdef}) shows that
each of $\alpha_d$, $\beta_d$, $\gamma_d$ and $\tau_d$ are in
${\Q}(\zeta_d)^+$. If $4$ does not divide $d$, then ${\Q}(\zeta_d)^+
= {\Q}(\zeta_{d'})^+$. Suppose now that $4 | d$, so $d' = d/2$ is
even. The Galois group
$\mathrm{Gal}({\Q}(\zeta_d)/{\Q}(\zeta_{d'}))$ is  generated by the
automorphism $\sigma$ which sends $\zeta_d$ to $\zeta_d^{1+d'} =
-\zeta_d$. In this case, each of the $p_j$ are odd since they are
relatively prime to $p$, $2|(d/2)$ and $(d/2)|p$.  It now follows
from (\ref{eq:greekdef}) that  each of $\alpha_d$, $\beta_d$ and
$\gamma_d$ are fixed by $\sigma$. Thus  these numbers are in  the
real subfield ${\Q}(\zeta_{d'})^+ = {\Q}(\zeta_{d'}) \cap
{\Q}(\zeta_d)^+$ in all cases.

We now prove the remaining assertions about $\alpha_d$. In
${\Q}(\zeta_d)$ we have
\begin{equation}
\label{eq:dummber}
\zeta_{d}^{p_j} - \zeta_{d}^{-p_j} = \zeta_d^{-p_j}(\zeta_d^{2 p_j} -1 )= \zeta_d^{-p_j}(
\zeta_{d'}^{p'_j} - 1)
\end{equation}
where  $d'= d$ and $p'_j = 2p_j$ if $d$ is odd, and  $d' = d/2$ and
$p'_j = p_j$  if $d$ is even. Recall that $p_j$ is prime to $p$, and
$d|2p$. Hence $p'_j = 2p_j$ is prime to $d' = d$ if $d$ is odd,
while $p'_j = p_j$ is prime to $d' = (d/2) | p$ if $d$ is even.
Thus $\zeta_{d'}^{p'_j}$ is a primitive ${d'}^{th}$ root of unity in
all cases.  In \cite[Prop. 2.8]{W} it is shown that $1 -
\zeta_{d'}^{p'_j}$ is a unit if $d'$ has at least two prime factors,
so $\alpha_d$ is a unit in this case.  Note that $d' > 1$ since $d >
2$.

The remaining possibility is that $d' = l^r$ for some prime $l$ and
some $r > 0$.  Then $1 - \zeta_{d'}^{p'_j}$ generates the unique
prime ideal $\mathcal{Q}$ over $l$ in ${\Q}(\zeta_{d'}) =
{\Q}(\zeta_{l^r})$ by \cite[p. 9]{W}.  From (\ref{eq:greekdef})  and
(\ref{eq:dummber}) we see that $\alpha_d \in {\Q}(\zeta_{d'})^+$ is
the product $\zeta \alpha'$ of a root of unity $\zeta \in
{\Q}(\zeta_d)$ with a generator $\alpha' \in {\Q}(\zeta_{d'})$ for
$\mathcal{Q}^4$.  Hence $\zeta \in {\Q}(\zeta_{d'})$, and $\alpha_d$
generates $\mathcal{Q}^4$ in ${\Q}(\zeta_{d'})$. The degree $e =
[{\Q}(\zeta_{d'}):{\Q}(\zeta_{d'})^+]$  is $2$ unless $d' = 2$, in
which case $e = 1$.  Since $\alpha_d \in {\Q}(\zeta_{d'})^+$, and
$\mathcal{P} \subset {\Q}(\zeta_{d'})^+$ generates the ideal
$\mathcal{Q}^e$ in ${\Q}(\zeta_{d'})$, we conclude that $\alpha_{d}$
generates the ideal $\mathcal{P}^2$ in ${\Q}(\zeta_{d'})^+$ unless
$d' = 2$, in which case $\alpha_{d}$ generates $\mathcal{P}^4 =
\mathcal{Q}^4$ in  ${\Q}(\zeta_{d'})^+ = {\Q}(\zeta_{d'}) = {\Q}$.
Since $d' = l^r > 1$, the extension
${\Q}(\zeta_{d})/{\Q}(\zeta_{d'})^+$ is totally ramified over $l$,
and has
 degree $2$ unless $d' = 2^r > 2$, in which case it has degree $4$.
\end{proof}

\begin{lemma}
\label{lem:nicest} With the notations of Lemma \ref{lem:dump}, each
of $\mathrm{Tr}_d(\frac{\beta_d}{\alpha_d})$, $\mathrm{Tr_d} \left
(\frac{1}{\alpha_d} \right )$, $\mathrm{Tr_d} \left
(\frac{\gamma_d}{\alpha_d} \right )$ and $\mathrm{Tr_d} \left
(\frac{\tau_d}{\alpha_d} \right )$
  lie in $2{\Z}$ unless $d' \in \{2,3, 5\}$. Moreover, for the
  remaining cases we obtain:
\begin{enumerate}
\item{} Suppose
$d' =2 $.  Then $\mathrm{Tr}_d(\frac{\beta_d}{\alpha_d}) = 0$,
$\mathrm{Tr}_d(\frac{1}{\alpha_d})$ generates the ideal
$\frac{1}{8}{\Z}$, and each of
$\mathrm{Tr}_d(\frac{\gamma_d}{\alpha_d}) = \frac{-4}{ \alpha_d}$
and $\mathrm{Tr}_d(\frac{\tau_d}{\alpha_d}) = \frac{4}{ \alpha_d}$
generate the ideal $\frac{1}{4}{\Z}$.
\item{}  Suppose
$d' =3 $.  Then $\mathrm{Tr}_d(\frac{\beta_d}{\alpha_d})$ lies in
$\frac{2}{9} {\Z}$, while $\mathrm{Tr}_d(\frac{1}{\alpha_d})$  and
$\mathrm{Tr}_d(\frac{\gamma_d}{\alpha_d}) =
\mathrm{Tr}_d(\frac{\tau_d}{\alpha_d})$  each generate the ideal
$\frac{2}{9}{\Z}$
\item{}Suppose $d' = 5$.  Then
$\mathrm{Tr}_d(\frac{1}{\alpha_d})$,
$\mathrm{Tr}_d(\frac{1}{\alpha_d})$,
$\mathrm{Tr}_d(\frac{\gamma_d}{\alpha_d})$ and
$\mathrm{Tr}_d(\frac{\tau_d}{\alpha_d})$ lie in $\frac{2}{5} {\Z}$.
\end{enumerate}
\end{lemma}

\begin{proof} Suppose first that $d' = 2$.  Then $d
= 4$, and each $p_j$ must be odd.  Now $\zeta_d^{p_j} +
\zeta_d^{-p_j} = i^{p_j} + i^{-p_j} = 0$, so $\beta_d = 0$ and
$\mathrm{Tr}_d(\frac{\beta_d}{\alpha_d}) = 0$. We have
\hbox{$[{\Q}(\zeta_d):{\Q}(\zeta_{d'})^+] = 2$}.  The number
$\alpha_d$ generates the ideal $2^4 {\Z}$ in ${\Q}(\zeta_{d'})^+ =
{\Q}$, $\gamma_d = \zeta_d^2 + \zeta_d^{-2} = -2$ and $\tau_d =
\zeta_d^4 + \zeta_d^{-4} = 2$.  Hence
$\mathrm{Tr}_d(\frac{1}{\alpha_d}) = \frac{2}{ \alpha_d}$ generates
the ideal $\frac{1}{8}{\Z}$, while
$\mathrm{Tr}_d(\frac{\gamma_d}{\alpha_d}) = \frac{-4}{ \alpha_d}$
and $\mathrm{Tr}_d(\frac{\tau_d}{\alpha_d}) = \frac{4}{ \alpha_d}$
each generate the ideal $\frac{1}{4}{\Z}$.

Now suppose $d' = 3$.  Then $d = 3$ or $d = 6$, so
\hbox{$[{\Q}(\zeta_d):{\Q}(\zeta_{d'})^+] = 2$}.  The number
$\alpha_d$ generates the ideal $3^2 {\Z}$ in ${\Q}(\zeta_{d'})^+ =
{\Q}$,and $\gamma_d = \zeta_d^2 + \zeta_d^{-2} = -1 = \zeta_d^4 +
\zeta_d^{-4} = \tau_d$. So $\mathrm{Tr}_d(\frac{\beta_d}{\alpha_d})
= \frac{2\beta_d}{ \alpha_d}$ lies in the ideal $\frac{2}{9}{\Z}$,
while $\mathrm{Tr}_d(\frac{1}{\alpha_d}) = \frac{2}{ \alpha_d}$ and
$\mathrm{Tr}_d(\frac{\gamma_d}{\alpha_d}) = \frac{-2}{ \alpha_d} =
\mathrm{Tr}_d(\frac{\tau_d}{\alpha_d})$  each generate the ideal
$\frac{2}{9}{\Z}$.

Next we consider $d' = 5$.  Then ${\Q}(\zeta_{d'})^+ =
{\Q}(\sqrt{5}) = L$ is quadratic over ${\Q}$, and $\alpha_d$
generates the square of the unique prime $\mathcal{P}$ over $l = 5$
in this extension.  Since $\mathcal{P}^2 $ is generated by $5$, we
conclude that for $\xi \in \{ \frac{\beta_d}{\alpha_d},
\frac{1}{\alpha_d}, \frac{\gamma_d}{\alpha_d},
\frac{\tau_d}{\alpha_d}\}$,  the number $5\xi$ is integral in $L$.
  This implies $\mathrm{Tr}_d(\xi) =
\mathrm{Tr}_{L/{\Q}}(2\xi) = \frac{2}{5}
\mathrm{Tr}_{L/{\Q}}(5\xi)\in\frac{2}{5} {\Z}$, as claimed.

 In the remaining computations we suppose $\xi \in
\{\beta_d,1,\gamma_d, \tau_d\}$, $d' > 3$ and $d' \ne 5$. If $d'$ is not a prime power
 then $\alpha_d$ is a unit, so
$\frac{\xi}{\alpha_d}$ is integral in ${\Q}(\zeta_{d'})^+$.  If $d'
= l^r$ for some
 prime $l$ and some $r > 0$, then either $l > 5$ or $r \ge 2$, and $\alpha_d$
 generates the square of the unique prime $\mathcal{P}$ over $l$ in
 ${\Q}(\zeta_{d'})^+$.
 In this case, the degree
of ${\Q}(\zeta_{d'})^+ = {\Q}(\zeta_{l^r})^+$ over ${\Q}$ is
$\phi(l^r)/2 = (l-1)l^{r-1}/2 \ge 3$,  so $\frac{\xi}{\alpha_d} \in
\mathcal{P}^{-2}$ and $\mathcal{P}^{-2}$ is contained in the inverse
different of  ${\Q}(\zeta_{d'})^+ = {\Q}(\zeta_{l^r})^+$ over
${\Q}$.  We conclude that in all cases,
 $\mathrm{Tr}_{d'+}(\frac{\xi}{\alpha_d}) \in {\Z}$, where $\mathrm{Tr}_{d'+}$
is the trace from ${\Q}(\zeta_{d'})^+$ to ${\Q}$. Hence
$$\mathrm{Tr}_d(\frac{\xi}{\alpha_d}) = \mathrm{Tr}_{d'+}(\frac{[{\Q}(\zeta_d):{\Q}(\zeta_{d'})^+]\xi}{\alpha_d}) =
[{\Q}(\zeta_d):{\Q}(\zeta_{d'})^+] \cdot
\mathrm{Tr}_{d'+}(\frac{\xi}{\alpha_d})$$  is in $2 {\Z}$ since
$[{\Q}(\zeta_d):{\Q}(\zeta_{d'})^+]$ is even.
\end{proof}

 Combining these results yields the following proposition, which proves
 Theorem \ref{thm:denomtheorem}.

\begin{proposition}
The rational numbers $T, R, S, $ and $U$ from Theorem
\ref{thm:denomtheorem} satisfy the following divisibility
properties:

\begin{enumerate}
\item If $(p,3)=(p,5)=1$, then $T, R, S,U$ are integers.
\item If $(p,3)=3$ and  $(p,5)=1$, then the denominators
are divisors of $9$.
\item If $(p,3)=1$ and $(p,5)=5$ , then the
denominators are  divisors of $5$.
\item If $(p,3)=3$ and $(p,5)=5$, then the denominators are
divisors of $45$.

\end{enumerate}
\end{proposition}

In terms of $T, R, S,U$ the Kreck-Stolz invariants of the lens
spaces are given by:

\begin{equation}\label{KST}
s_1=\frac{1}{2^5\cdot 7 \cdot p}(T+14S)\, ,\, s_2= \frac{1}{2^4\cdot
p}(R-S)\, ,\, s_3= \frac{1}{2^4\cdot p}(U-S)
\end{equation}

In order to determine the values of these Kreck-Stolz invariants on
a computer, we multiply $T, R, S,U$ by 45, find an integer
approximation, and use \ref{KST}.

\section{Examples}
\label{Examples}

Using a program written in Maple and C code we generate the
following lists of examples.  The program is available at www.math.upenn.edu/wziller/research ,  and
can be described briefly as follows.  For each given odd order
$r = |r(k,l)|=|\sigma_2(k)-\sigma_2(l)|<50000$ we produce a list of all
positively curved Eschenburg spaces with that given order of the
fourth cohomology group. To produce such a list, one needs to use four nested loops where the variable in each loop
 goes from 1 to 50000.
The list of such spaces is hence very large.
In the next step the program computes the
basic polynomial invariants $s(k,l)$ and $p_1(k,l)$ and  produces
a list of pairs whose basic invariants coincide. The program also
checks condition (C) and finds that it is always satisfied for
such pairs.
Generating the list and comparing the basic invariants are very time and memory intensive calculations which forced us to write the program in C code.
Surprisingly, there are only 437 pairs of spaces whose basic invariants coincide.

For this significantly smaller list of spaces the
Kreck-Stolz invariants are computed and compared, which can be done in Maple.

 We also indicate the cohomogeneity of the examples in the last column.
 Here $2+$ denotes the cohomogeneity two Eschenburg spaces with
$k_1 = k_2$, containing the $3$-Sasakian spaces,
 and $2-$ the case of $l_1 = l_2$.

We list the invariant $p_1 \in \Z_r$ as lying in the interval
$[0,r-1]$, the (orientation sensitive) invariants $s(k,l) \in
\Z_r$, which describes the linking form, as lying in
$(-\frac{r-1}{2} , \frac{r-1}{2}]$
 and  $s_1, s_2, s_{22} \in \Q/\Z$ as lying in $(-\frac{1}{2}, \frac12]$.  The
 advantage of choosing these intervals is that one sees
immediately when
the invariants just differ by a sign and hence the corresponding spaces are
orientation reversing homeomorphic or diffeomorphic.

We first produce a list of homotopy equivalent positively curved
Eschenburg spaces in  \taref{ET} for  $r\le 200$. Such examples
occur very frequently, e.g. there are $192$ such pairs for $r< 1000$.
See \cite{Sh} for the first examples of this type in the
literature.

 In order to find pairs of homeomorphic or diffeomorphic Eschenburg spaces we
  increased $r$ to $50000$.
There are $69$ homeomorphic pairs, the first 5 of which are listed in \taref{EH}, and
only  four diffeomorphic pairs, which we
 list in \taref{ED}. It is interesting to note  that for $r=26,973$ there are two
  Eschenburg metrics on the same manifold, one cohomogeneity two and the other
  cohomogeneity four. We do not know whether these two metrics are
  isometric or not, in particular whether the cohomogeneity four metric
  is in fact cohomogeneity two under a larger isometric  group action.

Next, we specialize to the subfamily of cohomogeneity two Eschenburg
spaces which are $3$-Sasakian, $E_{a,b,c}=\diag(z^a,z^b,z^c)
\backslash \SU(3) / \diag(z^{a+b+c},1,1)$. Homotopy equivalent
 $3$-Sasakian spaces again exist in abundance. We list the first 5
 examples with $r< 2000$ in \taref{ST}.

Our main goal for this subfamily was to produce two 3-Sasakian
spaces which are diffeomorphic to each other, see Theorem A. To do
this we increased $r$ to $10^7$.   Again we first find pairs of
3-Sasakian spaces where the basic invariants $s(k,l)$ and $p_1(k,l)$
coincide. The number of 3-Sasakian spaces with $r<10^7$ is  very
large, and this portion of the program was written in C code. The
program finds 3201  pairs with the same basic invariants, and we
then compute the Kreck-Stolz invariants.
  There are 96 pairs of homeomorphic $3$-Sasakian
spaces, the first five given in \taref{SH}, but
only one diffeomorphic pair, see \taref{SD}. A peculiar fact for
the homeomorphic pairs is that the sum $a+b+c$ and the orientation
is always the same,
which reduces the number  of possible pairs from 3201 to 152.
But we were not able to see this
directly.
 Notice also that this is not true for homotopy equivalent
pairs of $3$-Sasakian spaces, or for homeomorphic pairs of
Eschenburg spaces.


\providecommand{\bysame}{\leavevmode\hbox
to3em{\hrulefill}\thinspace}


\renewcommand{\arraystretch}{1.4}
\stepcounter{equation}
 \begin{table}[htbp]
   \begin{center}
      \begin{tabular}{|c|c|c|c|c|c|} \hline
 $r$ & $[k_1,k_2,k_3\;\vert\;l_1,l_2,l_3]$& $s$ & $s_{22}$ &
$p_1$  & Cohom \\
\hline \hline   43& $[21, 21, -2\;\vert\; 20, 20, 0]$ & 21& - 1/6 &26&1\\
\hline  43& $[8, 7, -5\;\vert\;6, 4, 0]$  & 21& -1/6  & 13 &4\\
\hline  &   &  &  &  & \\
\hline 101& $[50, 50, -2\;\vert\; 49, 49, 0]$ &  50 & -1/6 & 55 & 1\\
\hline 101& $[12, 10, -8\;\vert\; 9, 5, 0]$& 50 & -1/6 & 21 &4\\
\hline  &   &  &  & & \\
\hline 137& $[68, 68, -2\;\vert\; 67, 67, 0]$ & 68 & -1/6 & 73 & 1\\
\hline 137& $[19, 17, -7\;\vert\; 16, 13, 0]$ & 68 & -1/6 & 23 &4\\
\hline  &   &  &  & & \\
\hline 181& $[16, 16, -10\;\vert\; 13, 9, 0]$ & -26 & -1/6 & 85 & 2+ \\
\hline  181 & $[30, 26, -6\;\vert\;25, 25, 0]$ & 26 & 1/6 & 164 & 2- \\
\hline  &   &  &  & & \\
\hline 181& $[45, 43, -4\;\vert\; 42, 42, 0]$ & 43 & 0 & 89 & 2- \\
\hline 181& $[15, 14, -11\;\vert\; 12, 6, 0]$ & 43 & 0 & 35 &4\\
\hline
     \end{tabular}
   \end{center}
 \vspace{0.2cm}
     \caption{ Homotopy Equivalent Eschenburg Spaces for $r <
     200$}\label{ET}
      \end{table}


\stepcounter{equation}
\begin{table}[htpb]
\begin{center}
\begin{tabular}{|c|c|c|c|c|c|}
\hline
$r$ & $[k_1,k_2,k_3 \;\vert\; l_1,l_2,l_3]$& $s$ & $[p_1\,,\,s_2]$ &
$s_1$ & Cohom\\
\hline \hline   4001& $[79, 49, -50\;\vert\;  46, 32, 0]$& -1502 & [3336, -1043/8002]&49741/112028 & 4\\
\hline  4001&  $[75, 54, -51\;\vert\; 46, 32, 0]$ & 1502& [3336, 1043/8002] & 1877/8002 & 4\\
\hline  &   &  &  & & \\
\hline 8099 & $[71, 59, -94\;\vert\; 34, 2, 0]$ &  3085 &  [2184, -6975/32396] & -1055/9968 &4\\
\hline 8099 & $[92, 47,-85\;\vert\; 38, 16, 0]$ & -3085 &  [2184, 6975/32396] & -4285/9968 &4\\
\hline  &   &  &  & & \\
\hline 8671& $[83, 43, -96\;\vert\;  24, 6, 0]$ & 4216 &  [936, -11343/34684] & -941/10672 &4\\
\hline 8671&  $[97, 33, -88\;\vert\;  24, 18, 0]$ & - 4216 & [936, 11343/34684]  & -1417/74704 &4\\
\hline  &   &  &  & & \\
\hline 9889& $[104, 96, -86\;\vert\;  81, 33, 0]$ & 1719 & [65, 9505/39556] & 2961/79112 &4\\
\hline  9889 &$ [109, 101, -81\;\vert\;  81, 48, 0]$ & -1719 & [65, -9505/39556] & 275943/553784 &4\\
\hline  &   &  &  & & \\
\hline 11011&  $[144, 136, -76\;\vert\;  135, 69, 0]$ & -1899 & [5320, -6767/22022] & 31695/176176 &4\\
\hline 11011& $[152, 144, -68\;\vert\;  129, 99, 0]$ & -1899 &  [5320, -6767/22022]  & 12819/176176 &4\\
\hline
\end{tabular}
\end{center}
 \vspace{0.2cm}
     \caption{ Homeomorphic Eschenburg Spaces for $r <
     12000$}\label{EH}
 \end{table}

\stepcounter{equation}
 \begin{table}[htpb]
\begin{center}
\begin{tabular}{|c|c|c|c|c|c|} \hline
 $r$ & $[k_1,k_2,k_3 \;\vert\; l_1,l_2,l_3]$& $s$ &
$[p_1\,,\,s_2]$ &
$s_1$ & Cohom\\
\hline \hline   13361&  $[145, 121, -89\;\vert\; 113, 64, 0]$ & 1732 & [5905, 6839/53444] &-272959/748216&4\\
\hline  13361&   $[151, 127, -83\;\vert\; 104, 91, 0]$ & -1732 & [5905, -6839/53444] & 272959/748216 &4 \\
\hline  &   &  &  & & \\
\hline 26973 & $[154, 154, -158\;\vert\; 135, 15, 0]$ &  2119 & [5877, 123965/323676]& -6131/18648&2+\\
\hline 26973 &  $[389, 383, -67\;\vert\; 357, 348, 0]$ & -2119 &  [5877, -123965/323676]& 6131/18648&4\\
\hline  &   &  &  & & \\
\hline 35749&  $[185, 115, -186\;\vert\; 102, 12, 0]$ & 10989 &  [18648, 8920/35749] & -9018/35749&4\\
\hline 35749&   $[230, 111, -155\;\vert\; 108, 78, 0]$ & 10989 & [18648, 8920/35749] & -9018/35749&4\\
\hline  &   &  &  & & \\
\hline 42319& $[205, 141, -193\;\vert\;  114, 39, 0]$ & 7443 & [20142, 4123/84638] & -73317/677104 &4\\
\hline  42319 & $[191, 157, -195\;\vert\; 114, 39, 0]$ & -7443 & [20142, -4123/84638] & 73317/677104&4\\
\hline
\end{tabular}
\end{center}
 \vspace{0.2cm}
     \caption{ Diffeomorphic Eschenburg Spaces for $r\le
     50000$}\label{ED}
 \end{table}


\stepcounter{equation}
 \begin{table}[htpb]
\begin{center}
\begin{tabular}{|c|c|c|c|c|}
\hline
$r$ & $[a, b, c  \;\vert\;  a+b+c, 0, 0]$& $s$ & $s_{22}$ &
$p_1$ \\
\hline \hline   1267& $[316, 3, 1 \;\vert\; 320, 0, 0]$ & -319 & 1/3 & 813\\
\hline  1267& $[25, 19, 18 \;\vert\; 62, 0, 0]$ & -319 & 1/3  & 86 \\
\hline  &   &  &  &  \\
\hline 1277& $[181, 5, 2 \;\vert\; 188, 0, 0]$ &  533 & 1/6 & 453\\
\hline 1277& $[44, 19, 7 \;\vert\; 70, 0, 0]$ & -533 & -1/6 & 861\\
\hline  &   &  &  &  \\
\hline 1557& $[778, 1, 1 \;\vert\; 780, 0, 0]$ & 778 & 1/6 & 783\\
\hline 1557& $[139, 7, 4 \;\vert\;  150, 0, 0]$ & 778 & 1/6 & 1404\\
\hline  &   &  &  &  \\
\hline 1595& $[398, 3, 1 \;\vert\;  402, 0, 0]$ & -401 & 0 & 1018 \\
\hline  1595 & $ [36, 23, 13  \;\vert\; 72, 0, 0]$ & -401 & 0 & 798\\
\hline  &   &  &  &  \\
\hline 1619& $[105, 11, 4 \;\vert\; 120, 0, 0]$ & -237 & 0 & 1277\\
\hline 1619& $[132, 7, 5 \;\vert\; 144, 0, 0]$ & -237 & 0 & 997\\
\hline
\end{tabular}
\end{center}
\vspace{0.2cm}
     \caption{Homotopy equivalent 3-Sasakian Spaces $r<
     2000$}\label{ST}
 \end{table}

\stepcounter{equation}
 \begin{table}[htpb]
\begin{center}
\begin{tabular}{|c|c|c|c|c|cl}
\hline
$r$ & $[a,b,c \;\vert\; a+b+c,0,0]$& $s$ & $[p_1\,,\,s_2]$ &
$s_1$ \\
\hline \hline   28379 & $[171, 164, 1\;\vert\; 336, 0, 0]$ & -335 &[27139, -2393/56758] & -82869/3178448\\
\hline  28379&  $[223, 60, 53\;\vert\; 336, 0, 0]$ & -335 & [27139, -2393/56758]& -1104513/3178448 \\
\hline  &   &  &  &  \\
\hline 129503 & $[362, 291, 37\;\vert\; 690, 0, 0]$ &  12564 &  [45679, -80901/259006] & 69409/14504336\\
\hline 129503 &  $[423, 169, 98\;\vert\; 690, 0, 0]$& 12564 &  [45679, -80901/259006] & 5767541/14504336\\
\hline  &   &  &  &  \\
\hline 273581& $[717, 362, 13\;\vert\; 1092, 0, 0]$ & 91230 & [196280, 370663/1094324] & -393315/1094324\\
\hline 273581&   $[761, 241, 90\;\vert\; 1092, 0, 0]$ & 91230 & [196280, 370663/1094324] & 310179/1094324\\
\hline  &   &  &  &  \\
\hline 382025&  $[891, 368, 43\;\vert\; 1302, 0, 0]$ & -35741 &[334208, -294993/1528100] &-74669/436600\\
\hline  382025 &$[928, 191, 183\;\vert\; 1302, 0, 0]$ & -35741 & [334208, -294993/1528100] & 1442017/3056200\\
\hline  &   &  &  &  \\
\hline 442179&  $[1265, 347, 2\;\vert\; 1614, 0, 0]$ & -6448 & [346023, 115166/1326537] & -173889/611408\\
\hline 442179& $[1274, 311, 29\;\vert\; 1614, 0, 0]$ & -6448 & [346023, 115166/1326537]  & -21037/611408\\
\hline
\end{tabular}
\end{center}
\vspace{0.2cm}
     \caption{  Homeomorphic 3-Sasakian Spaces for $r<
     500000$}\label{SH}
 \end{table}

\vspace{6pt}

\stepcounter{equation}
 \begin{table}[htpb]
\begin{center}
\begin{tabular}{|c|c|c|c|c|}
\hline
$r$ & $[a,b,c \;\vert\; 0,0,a+b+c]$& $s$ & $[p_1\,,\,s_2]$ &
$s_1$ \\
\hline \hline   5143925 & $[2279, 1603, 384\;\vert\; 4266, 0, 0]$ & -1448517 & [390037, 36777/4115140]&-37291099/144029900\\
\hline   5143925 &  $[2528, 939, 799\;\vert\; 4266, 0, 0]$ & -1448517 &[390037, 36777/4115140] & -37291099/144029900 \\
\hline
\end{tabular}
\end{center}
\vspace{0.2cm}
     \caption{  Diffeomorphic 3-Sasakian Spaces for
     $r<10^7$}\label{SD}
 \end{table}

\end{document}